\documentclass{amsart}
\usepackage{amsmath,amssymb,latexsym,amsthm,amsfonts,amscd,epic,eepic,flafter,graphicx,verbatim, stmaryrd,mathrsfs,color}
\usepackage[all]{xy}
\author{Thomas Peters}
\address{Department of Mathematics, Columbia University, MC 4406\\ 2990 Broadway, New York,  NY 10027}
\email{tpeters@math.columbia.edu}

\title{\textbf{On L-spaces and non left--orderable 3-manifold groups}}
\begin{document}
\newtheorem{theorem}{Theorem}[section]
\newtheorem{lemma}[theorem]{Lemma}
\newtheorem{definition}[theorem]{Definition}
\newtheorem{proposition}[theorem]{Proposition}
\newtheorem*{thm*}{Theorem}
\newtheorem{rem}[theorem]{Remark}

\begin{abstract}
We show that a class of $3$--manifolds with non left--orderable fundamental group are Heegaard Floer homology $L$--spaces 
\end{abstract}

\maketitle

\section{Introduction}  
\label{section:intro}

Heegaard Floer homology is an powerful invariant of closed oriented 3-manifolds introduced by Ozsv\'ath and Szab\'o in \cite{OSzAnn1,OSzAnn2}.  It comes as a package of abelian groups, including a certain variant denoted, $\widehat{HF}(Y)$.  A $3$--manifold $Y$ is called an $L$--space if it is a rational homology $3$--sphere and its hat version of Heegaard Floer homology is ``as simple as possible" in the sense that the rank of $\widehat{HF}(Y)$ is equal to $|H_1(Y;\mathbb{Z})|$ (note that for a rational homology $3$--sphere, we always have $|H_1(Y;\mathbb{Z})|\leq rk\, \widehat{HF}(Y)$).  The class of $L$--spaces includes all lens spaces and is closed under connected sum as well as orientation reversal.  According to a theorem of N\'emethi, a three--manifold obtained as a plumbing of spheres is an $L$--space if and only if it is the link of a rational surface singularity \cite{nemethi}.  In fact, any $3$--manifold with spherical geometry is an $L$--space (see \cite{lspace}).  According to a theorem of Ozsv\'ath and Szab\'o, an $L$--space cannot have a co--orientable taut foliation \cite{OSzGenus}.  This provides a nice bridge between the world of pseudo-holomorphic curve invariants and the geometry of $3$-manifolds.  Though there is not yet a classification of $L$--spaces, there is a complete answer in the case of Seifert fibered spaces, according to a theorem of Lisca and Stipsicz \cite{LS} which states
\begin{theorem}\label{theorem:LS} Let $M$ be an oriented Seifert fibered rational homology $3$--sphere with base $S^2$.  Then the following statements are equivalent:
\begin{enumerate}
\item $M$ is an $L$--space
\item Either $M$ or $-M$ carries no positive, transverse contact structures
\item $M$ carries no transverse foliations
\item $M$ carries no taut foliations.
\end{enumerate}
\end{theorem}
Moreover, the existence of transverse foliations is completely understood and has a simple combinatorial answer given in terms of the Seifert invariants, by work of Eisenbud, Hirsch, Jankins, Neumann, and Naimi (see \cite{EHN}, \cite{JN}, \cite{JN2}, \cite{N}).

A group $G$ is called {\em left--orderable} if it may be given a strict total ordering $\prec$ which is left--invariant, i.e. $g\prec h$ if and only if $fg\prec fh$ if $f,g,h\in G$.  Non--left orderability has interesting consequences for the topology of a three--manifold.  For instance, Calegari and Dunfield showed that $3$--manifolds with non-left orderable fundamental group do not support co--orientable $\mathbb{R}$--covered foliations \cite{CD}.  Though in general, there is not a complete understanding of when a three--manifold group is left--orderable, Boyer, Rolfsen, and Wiest provide the answer in the case of Seifert fibered spaces \cite{BRW}:

\begin{theorem}\label{theorem:BRW} The fundamental group of a compact connected, Seifert fibered space $M$ is left--orderable if and only if $M\cong S^3$ or one of the following two sets of conditions holds:
\begin{enumerate}
\item $rank_\mathbb{Z}H_1(M)>0$ and $M\ncong\mathbb{P}^2\times S^1$;
\item $M$ is orientable, the base orbifold of $M$ is of the form $S^2(\alpha_1,\alpha_2,...,\alpha_n)$, $\pi_1(M)$ is infinite, and $M$ admits a transverse foliation.
\end{enumerate}
\end{theorem}

Putting together theorems \ref{theorem:LS} and \ref{theorem:BRW}, the remark that spherical manifolds are $L$--spaces, and the fact that closed Seifert fibered 3--manifolds with finite fundamental group are spherical, we see that a Seifert fibered $L$--space is either $S^3$ or its fundamental group is non--left orderable.   It is natural then, to explore this connection further.  Examples of infinite families of hyperbolic manifolds with non--left orderable fundamental group are provided by the work or Roberts, Shareshian, and Stein \cite{RSS}.  These manifolds were shown to be $L$--spaces by the work of Baldwin \cite{baldie2}.  In their paper \cite{CD}, Calegari and Dunfield determined that of the 128 closed hyperbolic manifolds of volume $<3$ which are $\mathbb{Z}/2$--homology spheres, at least 44 of them have non--left orderable fundamental group.  Dunfield later showed that all of these are in fact $L$--spaces \cite{dunfield}.  Further examples of non--left orderable three--manifold groups are provided by a paper of Dabkowski, Przytycki, and Togha \cite{DPT}. They prove

\begin{theorem}
\label{theorem:przy}
Let $\Sigma_n(L)$ denote the $n$--fold branched cyclic cover of the oriented link $L$, where $n>1$.  Then the fundamental group, $\pi_1(\Sigma_n(L))$, is not left-orderable in the following cases:
\begin{enumerate}
\item $L=T_{(2',2k)}$ is the torus link of type $(2,2k)$ with the anti-parallel orientation of strings, and $n$ is arbitrary. \\
\item $L=P(n_1,n_2,...,n_k)$ is the pretzel link of the type $(n_1,n_2,...,n_k)$, $k>2$, where either $n_1,n_2,...,n_k>0$ or $n_1=n_2=\cdots=n_{k-1}=2,n_k=-1$ and $k>3$.  The multiplicity of the covering is $n=2$.\\
\item $L=L_{[2k,2m]}$ is the $2$--bridge knot of type $p/q=2m+\frac{1}{2k}=[2k,2m]$, where $k,m>0$, and $n$ is arbitrary.\\
\item $L=L_{[n_1,1,n_3]}$ is the $2$--bridge knot of type $p/q=n_3+\frac{1}{1+\frac{1}{n_1}}$, where $n_1$ and $n_3$ are odd, positive integers.  The multiplicity of the covering is $n\leq3$.   
\end{enumerate}
 \end{theorem}

In this paper we show that

\begin{theorem}
\label{theorem:main}
All of the manifolds in theorem \ref{theorem:przy} are Heegaard Floer homology $L$--spaces.
\end{theorem}

Of course, the manifolds in (1) and (2), $\Sigma_3(4_1)$ from (4), and the covers of the trefoil from (4) are Seifert fibered and are hence covered by our previous remarks.  For the other cases, which are hyperbolic, we realize them as branched $double$ covers of of quasi--alternating links in the 3--sphere.  This allows us to apply a theorem of Ozsv\'ath and Szab\'o which states that the branched double cover of a quasi--alternating link in $S^3$ is an $L$--space.  We also give an independent proof for the manifolds in (1) also by realizing them as branched double covers of alternating pretzel links.

\subsection{Further Questions}
\label{section:questions}
We provide a list of unanswered questions which the author finds fascinating
\begin{enumerate}
\item We have a question of Ozsv\'ath and Szab\'o: is it true in general that a closed, oriented, and irreducible $3$-manifold is an $L$--space if and only if it has no co--orientable taut foliation?
\item What can one say in general about the relationship between Heegaard Floer homology and orderability properties of the fundamental group?
\item Given a knot or link, when is its $n$--fold cyclic cover an $L$--space?  For instance, it follows from Baldwin's classification of $L$--spaces among three--manifolds admiting genus one, one boundary component open books \cite{Baldie} that $\Sigma_n(3_1)$ is an $L$--space if and only if $n\leq5$ and $\Sigma_n(4_2)$ is an $L$--space for every $n$.
\item Which manifolds on the Hodgson--Weeks census are $L$--spaces?  Dunfield informs me that at least 3,000 of the 11,000 census manifolds are $L$--spaces \cite{dunfield}.
\item Is every $L$--space the branched double cover of a link in $S^3$? 
\item Give some description of hyperbolic $L$--spaces.
\item Connections to contact geometry:  Lisca and Stipsicz recently solved the existence problem for tight contact structues on Seifert fibered 3--manifolds \cite{LS2}:

\begin{theorem}
\label{theorem:LS2} A Seifert fibered 3--manifold admits a tight contact structure if and only if it is not orientation preserving diffeomorphic to the result of $(2n-1)$--surgery along the $(2,2n+1)$--torus knot $T_{2,2n+1}\subset S^3$ for some $n\in\mathbb{N}$.  
\end{theorem}

In proving this theorem, their classification of Seifert fibered $L$--spaces proved essential.  On the other hand, toroidal $3$--manifolds are known to admit infinitely many different contact structures (see \cite{CGH}, \cite{HKM}).  Little is known, however, about the existence of tight contact structures on hyperbolic $3$--manifolds.  Some information is provided by work of Baldwin \cite{baldie2}, and it is known that the Weeks manifold admits tight contact structures \cite{stipsicz}.  What can one say about tight contact structures on the manifolds from (3) and (4) in theorem \ref{theorem:przy}?

\end{enumerate}

\subsection{Acknowledgements.}
The author would like to thank his PhD supervisor, Peter Ozsv\'ath, for his continued guidance and support, and for suggesting the problem.  He would also like to thank Joshua Greene, Adam Levine, and Helge Moller Pederson for helpful conversations.  

\section{Background}
\label{section:review}

\subsection{A surgery presentation of the branched double cover of a link in $S^3$}
\label{section:surgery}
\begin{figure}[htbp] \begin{center}
\input{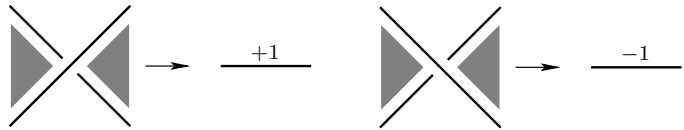}
\vspace*{-2mm} \caption{\label{fig:incidence} {\bf{Incidence assignment rules}}    }
\end{center}
\vspace*{-5mm}\end{figure}

\begin{figure}[htbp] \begin{center}
\input{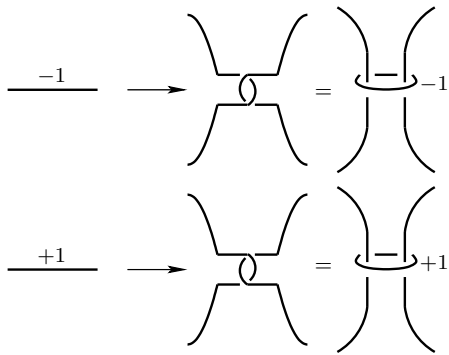}
\vspace*{-2mm} \caption{\label{fig:clasp} {\bf{Creating clasps out of incidences}}    }
\end{center}
\vspace*{-5mm}\end{figure}

\begin{figure}[htbp] \begin{center}
\input{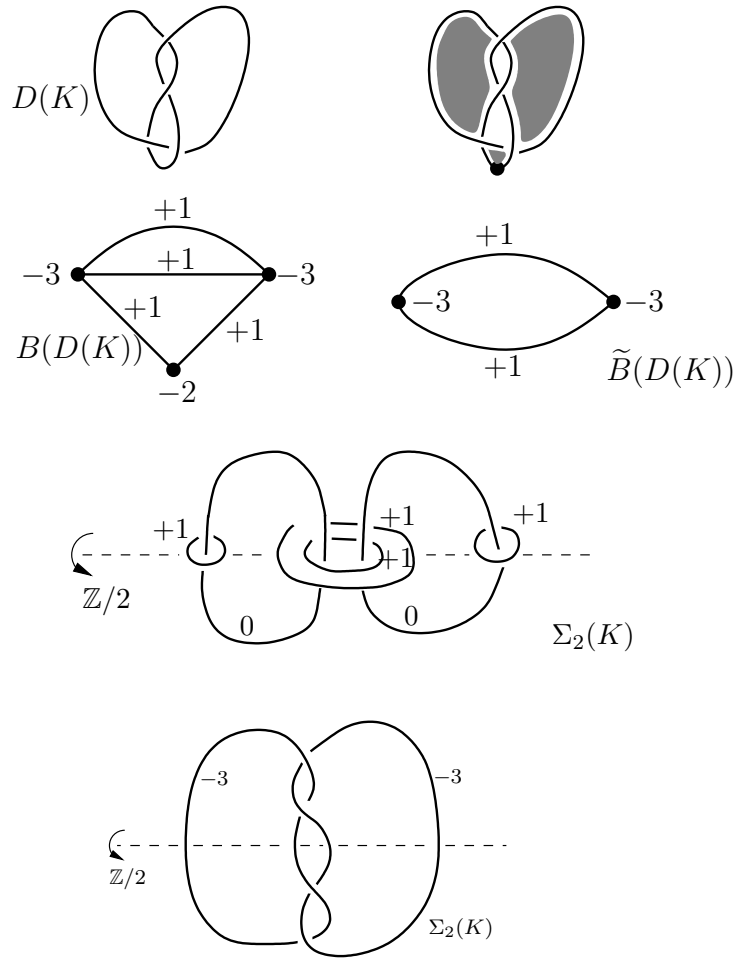}
\vspace*{-2mm} \caption{\label{fig:example} {\bf{Going from a diagram of the figure eight to a surgery presentation of its branched double cover}}    }
\end{center}
\vspace*{-5mm}\end{figure}

We begin with a review an algorithm that takes a diagram for a knot or link and produces a surgery presentation for its branched double cover, described in \cite{OSzGenus}.  Given a diagram of a link $D(K)$ pick an edge at random to mark.  Then checkerboard color the plane.  This allows us to produce the \emph{black graph} of $D(K)$, denoted $B(D(K))$: it is a planar graph whose vertices are in one-to-one correspondence with the black regions in our checkerboard coloring of the plane, and whose edges correspond to crossings in the diagram.  The edges are further decorated by an incidence number $\mu(e)=\pm1$ given by the rule of figure \ref{fig:incidence}.  The vertices are then weighted by the sum of the incidences of the incident edges $w(v) = -\sum_{e \text{ incident to } v}\mu(e)$.    We then form the \emph{reduced black graph}  $\widetilde{B}(D(K))$ by deleting the vertex which corresponds to the region touching the marked edge and then deleting all edges which are incident to this vertex.  We then draw a surgery diagram as follows:  for each vertex of $\tilde{B}$ we draw a planar unknot (such that all are unlinked).  For each edge between two vertices we add a right/left--handed clasp between the corresponding unknots according to the incidence of the edge (see figure \ref{fig:clasp}) or, equivalently, performing $\pm1$ surgery on an unknot which links the two components as shown in the figure.  If we chose to draw clasps, we frame each unknotted component by the weight on its corresponding vertex.  If we chose to draw linking $\pm1$ curves, then we $0$--frame each of the original unknots coming from the vertices of the reduced black graph.  This gives a surgery presentation for $\Sigma_2(K)$.  See figure \ref{fig:example} for an example.

Warning:  in this paper, we will occasionally work with decorated graphs such as the graph labeled $\widetilde{B}(D(K))$ in figure \ref{fig:example}.  Though asthetically similar, these diagrams are generally not the same as plumbing graphs (in the sense of \cite{Neumann}, for instance).  There is one exception, however:  when our graph is a tree and we delete the edge markings then we actually do have a plumbing description of our manifold.

One may visualize the involution on this manifold, giving rise to our link:  line up the $0$--framed circles on a line and then put in the $\pm1$ framed unknots each intersecting the axis of symmetry in $2$ points in such a way that the whole diagram has a symmetry about an axis, as shown in the example.  The complement of the surgered solid tori has the obvious branch locus--the axis drawn minus its intersections with the solid tori.  The involution of $S^3$ about this axis may be extended to a hyperelliptic involution of the surgery tori fixing longitudes and meridians setwise (but reversing their orientation).  It's easy to see that the quotient of the complement of the surgery solid tori under the involution is a ball minus a collection of disjoint sub-balls, one for each surgery torus (for instance, take as fundamental domain the ``upper half space" cut out by half solid tori).  This shows that the quotient orbifold is indeed topologically $S^3$.  We now determine its branch locus.  The branch locus in each of the solid tori is pair of arcs, each of which is isotopic (rel boundary) to a ``half" of the corresponding framing curve in the surgery diagram.  For instance, if we only had $n-1$ $0$--framed curves (no $\pm1$--surgeries), then the downstairs branch locus coming from the ``outside" (the complement of the solid tori) is a collection of $n$ arcs in the three--sphere.  Isotoping the downstairs branch loci (rel boundary) of the solid tori connects these arcs in such a way that we get a collection of $n$ unknots.  In a similar way, we can analyze what happens with the introduction of $\pm1$ surgeries, as above.  By pushing the branch loci into the ``outside" we see that $+1$ surgeries correspond to the introduction of right--handed ``crossings" between the corresponding unknots from above and that $-1$ surgeries correspond to the introduction of left--handed crossings.  Note that we're not quite working with knot diagrams and crossings--we're working with an unlink of some number of some components.  It still makes sense to talk about introducing ``crossings" between components once we chose a path between them: choose a ``parallel" copy of our path by displacing it slightly, keeping the endpoints on the link components.  With this framing curve in place, we can make sense of introducing ``crossings" between components.  In our case, the paths between components along which ``crossings" are added are the images of the cores of the surgery tori.  However, in our setting, we actually can interpret this as introducing crossings in a knot diagram (the plane can be taken to be the image of the half plane in the upper half space which intersects the axis of symmetry and the $0$--framed unknots union the boundary plane of the upper half space in the quotient...plus some stuff from the surgered solid tori).

\subsection{Quasi alternating knots.}
\label{section:qa}
In \cite{OSzDouble} Ozsv\'ath and Szab\'o defined the class of \emph{quasi--alternating} links--it is the smallest collection $\mathcal{Q}$ of links such that

\begin{figure}[htbp] \begin{center}
\input{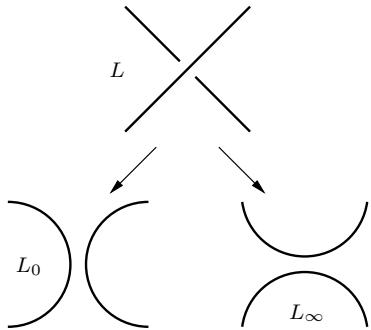}
\vspace*{-2mm} \caption{\label{fig:resolutions} {\bf{Resolving a crossing}}    }
\end{center}
\vspace*{-5mm}\end{figure}

\begin{itemize}
\item The unknot is in $\mathcal{Q}$.
\item If the link $L$ has a diagram with a crossing $c$ such that 
\begin{enumerate}
\item both resolutions of $c$, $L_0$ and $L_{\infty}$ as in figure \ref{fig:resolutions}, are in $\mathcal{Q}$,
\item $\det(L)=\det(L_0)+\det(L_{\infty})$,
\end{enumerate}
then $L$ is in $\mathcal{Q}$.
\end{itemize}

As in \cite{CK}, we shall call such a crossing $c$ as above a {\em quasi--alternating crossing} of
$L$ and that $L$ is {\em quasi--alternating at} $c$.

The class of quasi-alternating links extends the class of alternating links in the sense that if a link admits a connected alternating diagram then it is quasi--alternating.  In \cite{OSzDouble}, Ozsv\'ath and Szab\'o show that if a link $L$ is quasi--alternating, then its branched double cover is an $L$--space (though not every $L$--space arises in such a way, for instance $\Sigma(2,3,5)$ is an example).  

Quasi--alternating knots may be ``generated" by the following construction of Kofman and Champanerkar (see \cite{CK}).  Consider a crossing $c$ as above as a $2$--tangle with marked endpoints.  Let $\epsilon(c)=\pm1$ according to whether the overstrand has positive or negative slope.  We say that a rational $2$--tangle $\tau=C(a_1,...,a_m)$ {\em extends} $c$ if $\tau$ contains $c$ and $\epsilon(c)\cdot a_i\geq1$ for $i=1,...,m$.  They prove
\begin{theorem}
\label{theorem:twist}  Let $L$ be a quasi--alternating link with quasi--alternating crossing $c$ and $L'$ be obtained by replacing $c$ with an alternating rational tangle $\tau$ that extends $c$.  Then $L'$ is quasi--alternating at any crossing of $\tau$.
\end{theorem}

The reason that the branched double cover of a quasi--alternating link in $S^3$ is an $L$--spaces follows from the following construction of $L$--spaces (see \cite{lspace}).  Fix a closed, oriented three-manifold $Y$ and let $K$ be a framed knot in $Y$.  Then we have manifolds $Y_0$ and $Y_1$, obtained by $0$--surgery and $+1$--surgery on $K$.  We call the ordered triple $(Y,Y_0,Y_1)$ a {\em triad} of three--manifolds.  Suppose $Y,Y_0,Y_1$ are all rational homology three-spheres and $|H_1(Y)|=|H_1(Y_0)|+|H_1(Y_1)|$.  It follows from the surgery exact triangle in Heegaard Floer homology that if $Y_0$ and $Y_1$ are $L$--spaces, then so is $Y$.  The discussion in subsection \ref{section:surgery} shows that the branched double cover of a link and the branched double covers of its two resolutions at a crossing fit into a triad.  The previously mentioned theorem leads then to the inductive definition of quasi--alternating links.  A further consequence of the exact triangle shows that if $K\subset S^3$ is a knot in the three-sphere such that $S^3_r(K)$ is an $L$--space for some rational number $r>0$ (with respect to the Seifert framing) then $S^3_s(K)$ is an $L$--space for any rational $s>r$.    

\section{Proof of the theorem}
\label{section:proof}
\subsection{The manifolds in (1)}
\label{section:1}
Here we consider the manifolds $\Sigma_n(L)$ where $L=T_{(2',2k)}$ is the torus link of type $(2,2k)$ with the anti-parallel orientation of strings, and $n$ is arbitrary. 

\begin{figure}[htbp] \begin{center}
\input{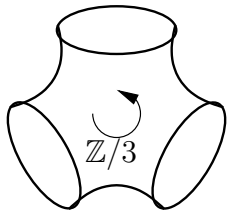}
\vspace*{-2mm} \caption{\label{fig:surface} {\bf{A $3$--fold branched covering of the disk branched over $2$ points downstairs.}}    }
\end{center}
\vspace*{-5mm}\end{figure}

Consider the genus $0$ open book decomposition of $S^3$ which has binding the braid axis for an unlink $\widetilde{L}$ of two components.  After performing $-\frac{1}{k}$--surgery on the binding of this open book, the unlink becomes $L$.  Each page of this open book meets $\widetilde{L}$ in exactly two points.  With this orientation, the $n$--fold strongly cyclic branched cover of the disk branched along two points is the $n$--times punctured sphere $S$.  The covering transormations consist of rotations through an axis which meets $S$ in two points through angles which are multiples of $2\pi/n$ and cyclically permute the boundary components (see figure \ref{fig:surface}).

\begin{figure}[htbp] \begin{center}
\input{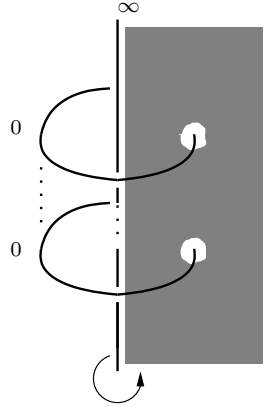}
\vspace*{-2mm} \caption{\label{fig:page} {\bf{An open book decomposition of $\#^{n-1}S^2\times S^1$.}}    }
\end{center}
\vspace*{-5mm}\end{figure}

\begin{figure}[htbp] \begin{center}
\input{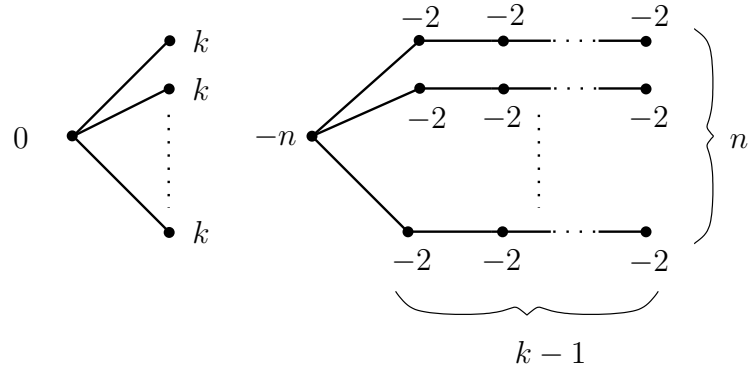}
\vspace*{-2mm} \caption{\label{fig:plumb1} {\bf{Two plumbing trees representing the manifolds from (1).  Both have $n$ branches.}}    }
\end{center}
\vspace*{-5mm}\end{figure}

\begin{figure}[htbp] \begin{center}
\input{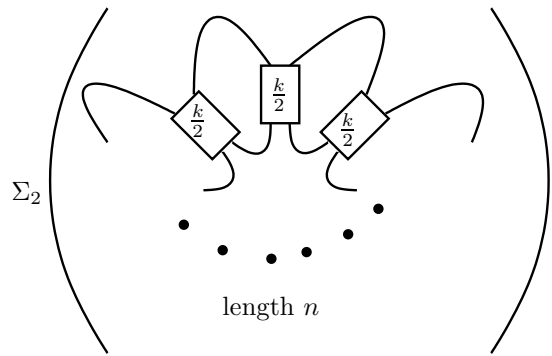}
\vspace*{-2mm} \caption{\label{fig:plumb2} {\bf{Realizing as a branched double cover.  Boxes represent full twists.}}    }
\end{center}
\vspace*{-5mm}\end{figure}

  The open book of $S^3$ with disk pages lifts in the $n$--fold branched cover to an open book with page $S$ and trivial monodromy--an open book decomposition of $\#^{n-1}S^2\times S^1$.  This open book decomposition is visualized in figure \ref{fig:page}.  Performing $-\frac{1}{k}$--surgery downstairs lifts to $-\frac{1}{k}$--surgery on the binding upstairs (with respect to the page framings).  This gives us the plumbing graph in the left hand side of figure \ref{fig:plumb1} (this is reached by beating one's head on figure \ref{fig:page}).  After a sequence of blowups/blowdowns, we reach the final plumbing graph, the right hand side of figure \ref{fig:plumb1} which, by the algorithm described in the background, may be  realized as a branched double cover of an alternating, hence quasi--alternating pretzel knot $P(\underbrace{k,k,...,k}_n)$ (see figure \ref{fig:plumb2}).  

\subsection{The manifolds in (3)}
\label{section:3}
We consider $\Sigma_n(L)$ for $L=L_{[2k,2m]}$ the $2$--bridge knot of type $p/q=2m+\frac{1}{2k}=[2k,2m]$, where $k,m>0$, and $n$ is arbitrary.
Surgery presentations of these manifolds may obtained by applying the ``Montesinos trick" or by appealing to a construction of Mulazzani and Vesnin \cite{MV}. 

\begin{figure}[htbp] \begin{center}
\input{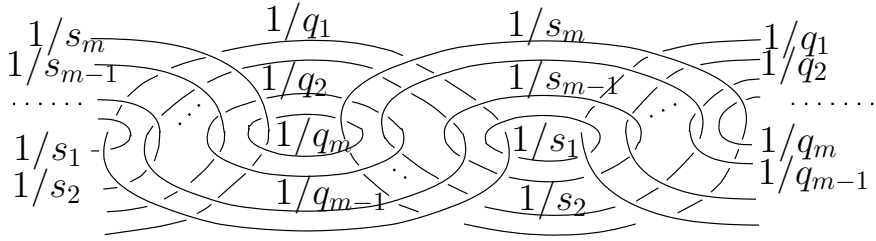}
\vspace*{-2mm} \caption{\label{fig:takahashi} {\bf{The $3$--manifold $T_{n,m}(1/q_j;1/s_j)$}}    }
\end{center}
\vspace*{-5mm}\end{figure}

Consider the $3$--manifold $T_{n,m}(1/q_j;1/s_j)$ defined by the surgery diagram in figure \ref{fig:takahashi}, with $2mn$ components, joined up to form a necklace.  Mulazzani and Vesnin prove that $T_{n,m}(1/q_j;1/s_j)$ is the $n$--fold cyclic branched covering of the two--bridge knot corresponding to the Conway parameters $[-2q_1,2s_1,...,-2q_m,2s_m]$.

\begin{figure}[htbp] \begin{center}
\input{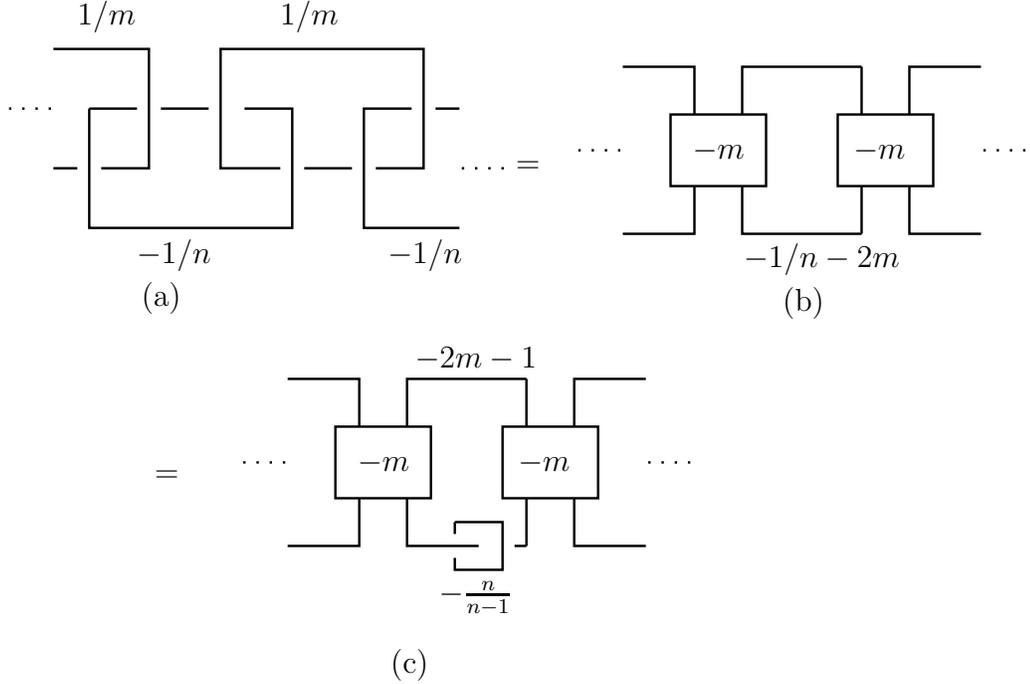}
\vspace*{-2mm} \caption{\label{fig:tak2} {\bf{Three views of the $3$--manifold $\Sigma_n(L_{[2k,2m]})$.  Each diagram is joined up to form a necklace of length $2k$.}}    }
\end{center}
\vspace*{-5mm}\end{figure}

\begin{figure}[htbp] \begin{center}
\input{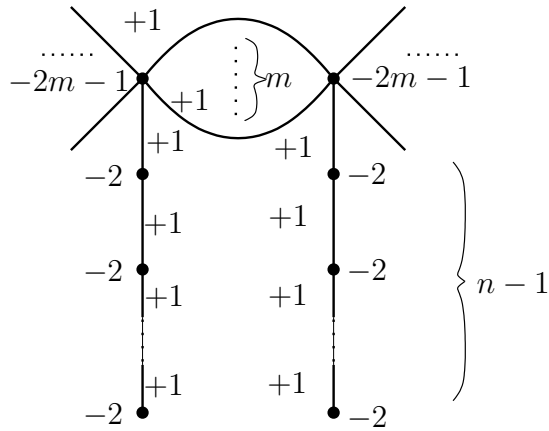}
\vspace*{-2mm} \caption{\label{fig:tak5}  {\bf{This is not a plumbing diagram!}}   }
\end{center}
\vspace*{-5mm}\end{figure}

Thus our manifolds have surgery presentation in figure \ref{fig:tak2} (a).  Twisting about each of the $\frac{1}{m}$--framed components give figures \ref{fig:tak2} (b) and \ref{fig:tak2} (c).  After some Kirby calculus we get figure \ref{fig:tak5} which we may realize as the branched double cover of the alternating knot shown in figure \ref{fig:tak6}.  

\begin{figure}[htbp] \begin{center}
\input{tak6.pstex_t}
\vspace*{-2mm} \caption{\label{fig:tak6}    }
\end{center}
\vspace*{-5mm}\end{figure}

\subsection{The manifolds in (4)}
\label{section:4}

\begin{figure}[htbp] \begin{center}
\input{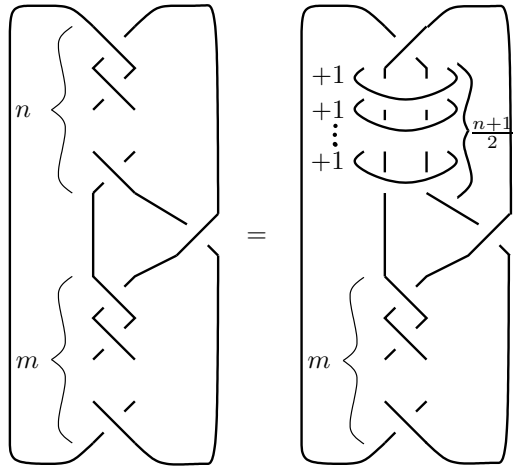}
\vspace*{-2mm} \caption{\label{fig:Lnm} {\bf{The knots $L_{[n,1,m]}$}}    }
\end{center}
\vspace*{-5mm}\end{figure}

\begin{figure}[htbp] \begin{center}
\input{Lnm2.pstex_t}
\vspace*{-2mm} \caption{\label{fig:Lnm2} {\bf{The knots $L_{[n,1,m]}$}}    }
\end{center}
\vspace*{-5mm}\end{figure}

Finally, we consider the manifolds $\Sigma_k(L)$ where $L=L_{[n,1,m]}$ is the $2$--bridge knot of type $p/q=m+\frac{1}{1+\frac{1}{n}}$, where $n$ and $m$ are odd, positive integers.  The multiplicity of the covering is $n\leq3$.  For $n=2$, these manifolds are lens spaces, and hence $L$--spaces, so we consider the case of $n=3$.  To construct surgery diagrams for these manifolds, we could appeal to the work of Mulazzani and Vesnin.  Instead, we use the ``Montesinos trick" (see \cite{Rolfsen}), as in figures \ref{fig:Lnm}, \ref{fig:Lnm2}, and \ref{fig:Lnm3}. 

\begin{figure}[htbp] \begin{center}
\input{Lnm3.pstex_t}
\vspace*{-2mm} \caption{\label{fig:Lnm3} {\bf{The manifold $\Sigma_3(L_{[n,1,m]})$}.}    }
\end{center}
\vspace*{-5mm}\end{figure}

\begin{figure}[htbp] \begin{center}
\input{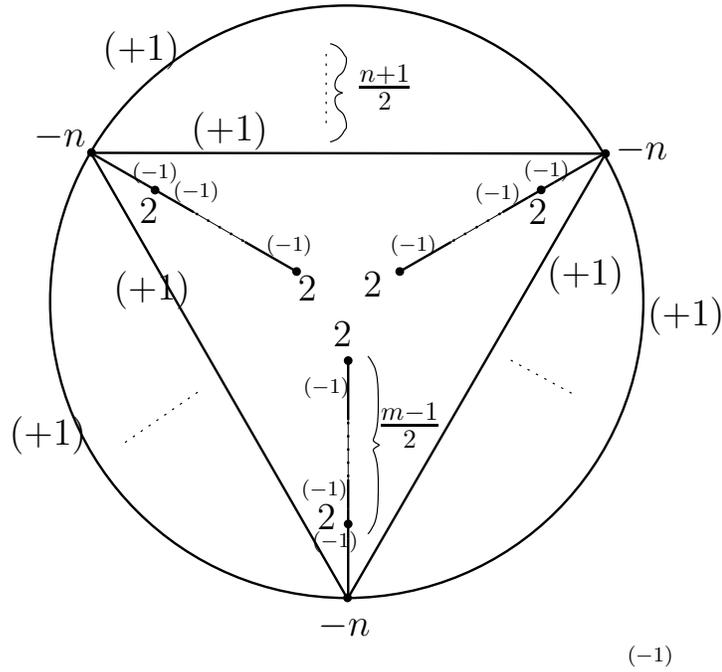}
\vspace*{-2mm} \caption{\label{fig:Lnm5} {\bf{The manifold $\Sigma_3(L_{[n,1,m]})$}}    }
\end{center}
\vspace*{-5mm}\end{figure}

\begin{figure}[htbp] \begin{center}
\input{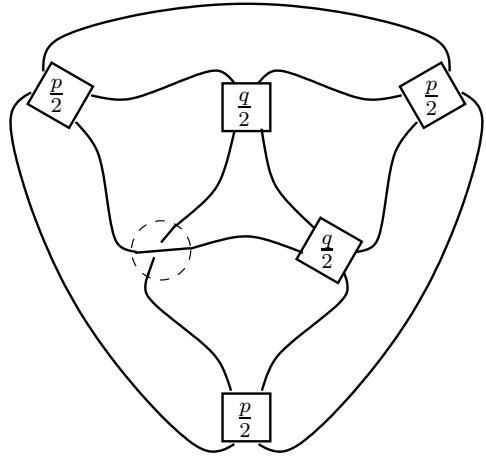}
\vspace*{-2mm} \caption{\label{fig:L} {\bf{The link $L$}}    }
\end{center}
\vspace*{-5mm}\end{figure}

 After some blowups/blowdowns, we see that our surgery diagram can be represented by the graph in figure \ref{fig:Lnm5}, which we claim is (usually) the reduced black graph of a quasi--alternating knot.  Indeed, consider the family of links $K(p,q)$ ($p,q\in \mathbb{N}_{\geq0}$) shown below.  Then our most recent diagram is the branched double cover of the link $K(\frac{n+1}{2},\frac{m+1}{2})$ (add a ghost vertex at the center, connected to the closest $2$'s by 3 $-1$--marked edges).  $K(1,1)$ is the $3$--braid $(\sigma_1\sigma_2)^3$ which, though not quasi--alternating, has an $L$--space as branched double cover (see \cite{Baldie}).  This manifold can be realized as the triple branched cover of the trefoil knot, which is the spherical space form $\mathbb{S}^3/Q_8$.  We claim that $K(p,q)$ quasi--alternating if $pq>1$.  By a symmetry of the diagram for $K(p,q)$, we may assume that $p>1$.  Using theorem \ref{theorem:twist}, it is enough to see that the link $L$ (shown in figure \ref{fig:L}) is quasi--alternating at the circled crossing.  Figures \ref{fig:0res} and \ref{fig:infres} show that the two resolutions $L_0$ and $L_\infty$ admit connected alternating diagrams and hence are quasi--alternating.  

\begin{figure}[htbp] \begin{center}
\input{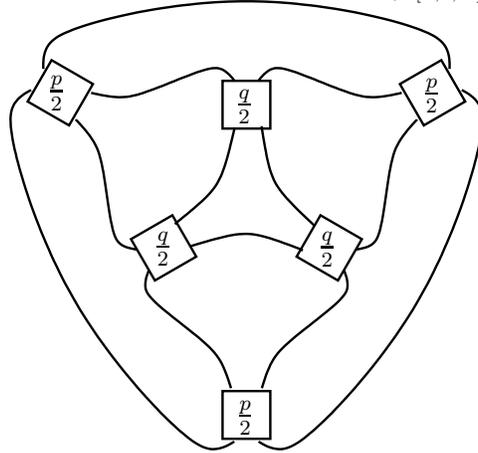}
\vspace*{-2mm} \caption{\label{fig:pqrstu} {\bf{The family of links $K(p,q)$}}    }
\end{center}
\vspace*{-5mm}\end{figure}

\begin{figure}[htbp] \begin{center}
\input{0res.pstex_t}
\vspace*{-2mm} \caption{\label{fig:0res} {\bf{$L_0$}}    }
\end{center}
\vspace*{-5mm}\end{figure}

\begin{figure}[htbp] \begin{center}
\input{infres.pstex_t}
\vspace*{-2mm} \caption{\label{fig:infres} {\bf{$L_\infty$}}    }
\end{center}
\vspace*{-5mm}\end{figure}

Let \[
	 A= \begin{vmatrix} 
1-2p & p & p &  &  &  &  &  &  &  &\\
p & 1-2p & p & -1 &  &  &  &  &  &  & \\
p & p & 1-2p &  &  & \cdots &  & -1 &  &  & \\
 & -1 &  & 2 & -1 &  &  &&  &  & \\
 &  &  & -1 & 2 &  &  & &  &  & \\
&&\vdots&&&\ddots&-1&&& &\\
&&&&&-1&2&&&\\
 &  & -1 &  &  &  &  & 2& -1 & &\\
&&&&&&&-1&2&&\\
&&&&&&&&&\ddots&-1 \\
&&&&&&&&&-1&2
\end{vmatrix}\]

Where both strings of $2$'s have length $q-1$.  Then $|A| = \det L$.

Now let \[
	 B= \begin{vmatrix} 
-2p & p & p &  &  &  &  &  &  &  &\\
 p& 1-2p & p &  &  & \cdots &  & -1 &  &  & \\
 p& p & 1-2p & -1 &  &  &  &  &  &  & \\
 &  & -1 & 2 & -1 &  &  &&  &  & \\
 &  &  & -1 & 2 &  &  & &  &  & \\
\vdots&&&&&\ddots&-1&&& &\\
&&&&&-1&2&&&\\
-1 &  &  &  &  &  &  & 2& -1 & &\\
&&&&&&&-1&2&&\\
&&&&&&&&&\ddots&-1 \\
&&&&&&&&&-1&2
\end{vmatrix}\]

Where again both strings of $2$'s have length $q-1$.  Then $|C| = \det L_0$.  Finally let\[
	 C = \begin{vmatrix} 
 1-2p & p &  &  & \cdots &  & -1  &  &  & \\
 p & 1-2p & -1 &  &  &  &  &  &  & \\
  & -1 & 2 & -1 &  &  &&  &  & \\
  &  & -1 & 2 &  &  & &  &  & \\
&\vdots&&&\ddots&-1&&& &\\
&&&&-1&2&&&\\
   & -1 &  &  &  &  & 2& -1 & &\\
&&&&&&-1&2&&\\
&&&&&&&&\ddots&-1 \\
&&&&&&&&-1&2
\end{vmatrix}\]  Where again the first string of $2$'s is length $s-1$ and the second string is of length $t-1$.  Then $|C| = \det L_\infty$.  Clearly we have $A = B+C$.  We claim that $B,C>0$ so that $|A| = |B|+|C|$ and $L$ is quasi-alternating.  Let's start with $B$:  consider the more general matrix
 \[
	 B(p,q,r)= \begin{vmatrix} 
-2p & p & p &  &  &  &  &  &  &  &\\
 p& 1-2p & p &  &  & \cdots &  & -1 &  &  & \\
 p& p & 1-2p & -1 &  &  &  &  &  &  & \\
 &  & -1 & 2 & -1 &  &  &&  &  & \\
 &  &  & -1 & 2 &  &  & &  &  & \\
&\vdots&&&&\ddots&-1&&& &\\
&&&&&-1&2&&&\\
 & -1 &  &  &  &  &  & 2& -1 & &\\
&&&&&&&-1&2&&\\
&&&&&&&&&\ddots&-1 \\
&&&&&&&&&-1&2
\end{vmatrix}\]
Where the first string of $2$'s is of length $q$ and the second of length $r$.  Then $B = B(p,q-1,q-1)$.  We claim that $b(p,q,r):=\det(B(p,q,r))>0$ for any $p>0$.  A simple calculation shows that $b(p,q,r)$ satisfies the recurrences
\[
	b(p,q,r) = 2b(p,q-1,r)-b(p,q-2,r), q>1
\]
\[
	b(p,q,r) = 2b(p,q,r-1)-b(p,q,r-2), r>1
\]
which leads to the solution
\begin{align*}
	b(p,q,r) = &rq(b(p,1,1)-b(p,0,1)-b(p,1,0)+b(p,0,0))\\
	&+r(b(p,0,1)-b(p,0,0))+q(b(p,1,0)-b(p,0,0)) +b(p,0,0).
\end{align*}
A little computer assistance shows then that
\[
	b(p,q,r) = rq(0)+r(3p^2)+q(3p^2)+(-2p+6p^2)
\]
which is positive for any $p>0$.  Similarly consider the matrices
\[	 C(p,q,r) = \begin{vmatrix} 
 1-2p & p &  &  & \cdots &  & -1  &  &  & \\
 p & 1-2p & -1 &  &  &  &  &  &  & \\
  & -1 & 2 & -1 &  &  &&  &  & \\
  &  & -1 & 2 &  &  & &  &  & \\
\vdots&&&&\ddots&-1&&& &\\
&&&&-1&2&&&\\
   -1&  &  &  &  &  & 2& -1 & &\\
&&&&&&-1&2&&\\
&&&&&&&&\ddots&-1 \\
&&&&&&&&-1&2
\end{vmatrix}\]
Where the first string of $2$'s is length $q$ and the second of length $r$.  Then $C = C(p,q-1,q-1)$.  We claim that $c(p,q,r):=\det(C(p,q,r))>0$ for any $p>1$.  Similar to before, we see that $c(p,q,r)$ satisfies the recurrences
\[
	c(p,q,r) = 2c(p,q-1,r)-c(p,q-2,r), q>1
\]
\[
	c(p,q,r) = 2c(p,q,r-1)-c(p,q,r-2), r>1
\]
which leads to the solution
\begin{align*}
	c(p,q,r) = &rq(c(p,1,1)-c(p,0,1)-c(p,1,0)+c(p,0,0))\\
	&+r(c(p,0,1)-c(p,0,0))+q(c(p,1,0)-c(p,0,0)) +c(p,0,0).
\end{align*}
Hence
\[
	c(p,q,r)=rq(3p^2)+r(-2p+3p^2)+q(-2p+3p^2)+(1-4p+3p^2)
\]
which is clearly positive for $p>1$.

\end{document}